\def\draft{n}
\documentclass{amsart}
\usepackage{fullpage,amssymb,epic,eepic,epsfig,amscd,pb-diagram,lamsarrow,pb-lams}

\theoremstyle{plain}

\newtheorem{theorem}{Theorem}
\newtheorem{proposition}{Proposition}[section]
\newtheorem{lemma}[proposition]{Lemma}

\theoremstyle{definition}
\newtheorem{definition}[proposition]{Definition}

\theoremstyle{remark}

\newtheorem{remark}[proposition]{Remark}

\def\printname#1{
	\if\draft y
		\smash{\makebox[0pt]{\hspace{-0.5in}
			\raisebox{8pt}{\tt\tiny #1}}}
	\fi
}

\newcommand{\psdraw}[2]
         {\begin{array}{c} \hspace{-1.3mm}
	\raisebox{-4pt}{\epsfig{figure=draws/#1.eps,width=#2}}
	\hspace{-1.9mm}\end{array}}

\newlength{\standardunitlength}
\catcode`\@=11
\long\def\@makecaption#1#2{%
     \vskip 10pt

\setbox\@tempboxa\hbox{
       \small\sf{\bfcaptionfont #1. }\ignorespaces #2}%
     \ifdim \wd\@tempboxa >\captionwidth {%
         \rightskip=\@captionmargin\leftskip=\@captionmargin
         \unhbox\@tempboxa\par}%
       \else
         \hbox to\hsize{\hfil\box\@tempboxa\hfil}%
     \fi}
\font\bfcaptionfont=cmssbx10 scaled \magstephalf
\newdimen\@captionmargin\@captionmargin=2\parindent
\newdimen\captionwidth\captionwidth=\hsize
\catcode`\@=12

\newcommand{\tr}{\operatorname{tr}}

\def\lbl#1{\label{#1}\printname{#1}}

\def\BZ{\mathbb Z}
\def\BQ{\mathbb Q}

\def\BC{\mathbb C}

\def\B{\mathcal B}

\def\D{\Delta}

\def\O{\mathcal O}

\def\N{\mathcal N}

\def\R{\mathcal R}

\def\La{\Lambda}

\def\S{\Sigma}

\def\ga{\gamma}

\def\ihs{integral homology 3-sphere}

\def\la{\langle}
\def\ra{\rangle}

\def\logp{\log'}      
\def\cyclic{\mathrm{cyclic}}

\def\e{\epsilon}

\def\lk{\mathrm{lk}}

\def\pt{\partial}

\def\sminus{\smallsetminus}

\def\ti{\widetilde}

\def\mtwo#1#2#3#4{ \begin{pmatrix}
                   #1 & #2 \\
                   #3 & #4 
                   \end{pmatrix} }

\def\Lxx{L^{xx}}
\def\Lxy{L^{xy}}
\def\Lyx{L^{yx}}
\def\Lyy{L^{yy}}
\def\fx{\bar{x}}
\def\ft{\bar{t}}

\def\chiD{\chi_{\Delta}}
\def\ygraph{$\mathrm{Y}$-graph}

\def\Mat{\mathrm{Mat}}

\def\Lhat{\hat\Lambda}
\def\Lloc{\Lambda_{\mathrm{loc}}}

\def\Zrat{Z^{\mathrm{rat}}}

\def\links{\mathrm{links}}
\def\NO{\N(\O)}

\def\Bla{\B(\La\to\BZ)}

\def\longto{\longrightarrow}

\def\pt{\partial}

\begin{document}


\title[A surgery view of boundary links]{A surgery view of boundary links}

\author{Stavros Garoufalidis}
\address{School of Mathematics \\
          Georgia Institute of Technology \\
          Atlanta, GA 30332-0160, USA. }
\email{stavros@math.gatech.edu}
\author{Andrew Kricker}
\address{Department of Mathematics \\ 
         University of Toronto \\
         Toronto, Ontario, Canada M5S 3G3.}
\email{akricker@math.toronto.edu}

\thanks{S.G. was partially supported by an NSF grant
        DMS-02-03129 and by an Israel-US BSF grant. A.K. was supported by
        an Israel Fellowship. \newline
        This and related preprints can also be obtained at
{\tt http://www.math.gatech.edu/$\sim$stavros }. 
\newline
1991 {\em Mathematics Classification.} Primary 57N10. Secondary 57M25.
}

\date{
This edition: October 16, 2002 \hspace{0.5cm} First edition: April 29, 2002.}


\begin{abstract}
A celebrated theorem of Kirby identifies the set of closed oriented connected
3-manifolds with the set of framed links in $S^3$ modulo two moves.
We give a similar description for the set of knots (and more
generally, boundary links) in homology 3-spheres. As an application, we define
a noncommutative version of the Alexander polynomial of a boundary link.
Our surgery view of boundary links is a key ingredient in a construction
of a rational version of the Kontsevich integral, which is described
in subsequent work.
\end{abstract}

\maketitle



\section{Introduction}
\lbl{sec.intro}

Surgery (or cut-and-paste topology) is a method of modifying a manifold
to another one. Surgery was successfully used in the sixties to geometrically
realize algebraic invariants of manifolds and leads, for example, to
a classification of high dimensional manifolds of a fixed homotopy type.

Surgery has also been fruitfully applied to the case of embedding questions, 
most
notably to codimension 2 embeddings, i.e., knot theory. For an excellent
survey, see \cite{LO}. This
was pioneered by Levine, who used surgery to geometrically realize known
knot (and link) invariants, such as the Alexander polynomial, the Alexander
module, and well-known concordance invariants of knots, \cite{Le1}.
Rolfsen used surgery for similar reasons in his reader-friendly introduction
to knot theory, \cite{Rf}. 
The key idea behind this is the fact that knots (or rather, knot projections)
can be unknotted via a sequence of crossing changes, and that a crossing
change can be achieved by surgery on a $\pm 1$-framed unknot as follows:

\begin{figure}[htpb]
$$ 
\psdraw{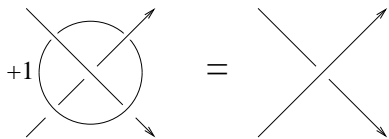}{2in} 
$$
\caption{A crossing change can be achieved by surgery on a unit framed
unknot.}\lbl{split}
\end{figure}

Thus, every knot $K$ in $S^3$ can be obtained by surgery on a framed link $C$
in the complement $S^3\sminus\O$ of an unknot $\O$. We will call such 
a link $C$, an {\em untying link} for $K$. Observe that untying links are
framed, and null homotopic in $S^3\sminus\O$, the interior of a solid torus.
Further, their linking matrix is invertible over $\BZ$, since surgery on
them gives rise to a \ihs .

In the above mentioned applications of surgery theory to knot
theory, one starts with a knot invariant (such as the Alexander polynomial) 
with a known behavior under surgery. Using this information, one can realize
all values of the invariant by choosing suitable untying links on 
$S^3\sminus\O$ and performing surgery on them. This is exactly what Levine, 
Rolfsen and others did.

There is one point, though, that the above literature does not usually
discuss: a 
knot can be untied by many different links in $S^3$. This is due 
to the fact that surgery is related to handlebodies that have different
descriptions. If two framed links are related by a $\kappa$ ($=\kappa_1$ or
$\kappa_2$) move, then surgery on them gives rise to diffeomorphic
manifolds. Here, $\kappa_1$ corresponds to adding a parallel of a link 
component to another, and $\kappa_2$  corresponds to adding/removing a 
unit-framed unlink away from a link.

The above discussion applies in the case of knots $K$ in {\em \ihs s}
$M$ (that is, closed 3-manifolds with the same $\BZ$-homology as $S^3$), and
leads to the following {\em surgery map}:

\begin{equation}
\lbl{eq.surgery}
\NO/\la\kappa\ra\longrightarrow \text{Knots},
\end{equation}
where $\mathrm{Knots}$ is the set of knots in homology spheres, and 
$\NO$ denotes the set of null homotopic links $C$ with $\BZ$-invertible 
linking matrix in the complement of a standard unknot $\O$ in $S^3$. 
The next theorem explains what we mean by the {\em surgery view of knots}.

\begin{theorem}
\lbl{thm.surgeryK}
The surgery map \eqref{eq.surgery} is 1-1 and onto.
\end{theorem}

The technique of proving this theorem is well-known in the seventies;
in a sense it is a relative version of Kirby's theorem for manifolds with
boundary. The reader may ask why did the classical literature not 
state the above result, and why is it of any use nowadays? As we mentioned 
already, if one is interested in {\em realizing} values of known knot 
invariants, one needs to know that the surgery map is onto. If one is 
interested in {\em constructing} a knot invariant, starting from an invariant 
of null homotopic
links in the complement of an unknot, then one needs to know the kernel
of the surgery map. This is precisely the situation that we are interested in.
In forthcoming work \cite{GK} we will use the above theorem in 
order to construct a new knot invariant, namely a rational version $\Zrat$ of 
the Kontsevich integral of a knot. This work motivated us to formulate
Theorem \ref{thm.surgeryK}.

In a sense, Theorem \ref{thm.surgeryK} defines knot theory in terms of link
theory in a solid torus. One might think that this is a step in the wrong 
direction. However, this view of knots is closely related to the {\em Homology
Surgery} view of knots of Cappell-Shaneson, \cite{CS1,CS2}. Cappell and 
Shaneson
study the closed 3-manifold obtained by $0$ surgery on a knot. This manifold
has the homology of a circle but not the $\BZ$-equivariant homology of a 
circle.

In \cite{GK}, we actually construct a rational invariant of $F$-links.
Boundary links ($\pt$-links), and their cousins, $F$-links are a 
generalization of knots which we now explain. Recall first that
every knot $K$ in a homology sphere $M$ gives rise to a map 
$\phi:\pi_1(M\sminus K)\to H_1(M\sminus K)\cong\BZ$ which maps a meridian
to a generator of $\BZ$. An $F$-{\em link} is a triple $(M,L,\phi)$ of
a link $L$ in a homology sphere $M$, where $\phi:\pi_1(M\sminus L)\to F$
is a map that sends meridians of $L$ to generators of the free group $F$.
Here, $F$ denotes the free group in generators $t_1,\dots,t_g$, where $r$
is the number of components of $L$. An $F$-link gives rise to a map
$M\sminus L \to \vee^r S^1=K(F,1)$ and the preimage of generic points
on each circle gives rise to a disjoint union of surfaces $\S_i$ with
boundary the $i$th component $L_i$ of $L$. Thus, $(M,L)$ is a {\em boundary
link}. 

It turns out that every boundary link has an $F$-structure, however not
a canonical one. For a discussion of this, we refer the reader to \cite{GL2}
and \cite{GK}.

Generalizing Theorem \ref{thm.surgeryK}, let $\O$ denote a standard unlink
of $r$ components in $S^3$, with a canonical isomorphism $\phi:
\pi_1(S^3\sminus\O) \cong F$. If $C$ is a null homotopic link in 
$S^3\sminus \O$
with $\BZ$ invertible linking matrix, then after surgery it gives rise to
an $F$-link $(M,L,\phi)$. The next theorem describes the {\em surgery view
of $F$-links}.

\begin{theorem}
\lbl{thm.surgeryF}
The surgery map gives a 1-1 onto correspondence 
$$
\NO/\la\kappa\ra\longleftrightarrow F-\links .
$$
\end{theorem}

As we mentioned above, we will use this theorem in \cite{GK} in 
order to construct a rational version of the Kontsevich integral. At present, 
we will use the above theorem to define a noncommutative version of the 
Alexander polynomial of an $F$-link, as well as a matrix-valued invariant of 
$F$-links which is essentially equivalent to the Blanchfield pairing of their 
free cover.

Consider a link $C \in \NO$. Since $C \subset X_0=S^3\sminus\O$ is null 
homotopic, it lifts to a link $\ti C$ in the free $F$-cover $\ti X_0$ of $X_0$.
Since $\ti X_0$ is contractible, it follows that we can define linking numbers
between the components of $\ti C$, and thus we can define the {\em equivariant
linking matrix} $\lk(\ti C)$ of $C$ by:
$$
\lk(\ti C)_{ij}=\left( \sum_{g \in F} \lk(\ti C_i, g \, 
\ti C_j) \, g \right)_{ij}
$$
where $\ti C_i$ denotes a lift of a component $C_i$ of $C$. The above sums
are finite and $\lk(\ti C)$ is a matrix over the group-ring $\La=\BZ[F]$
of $F$. The size of the matrix is the number of components of $C$.
The equivariant linking matrix of $C$ is a {\em Hermitian
matrix} (that is, a matrix $A$ over $\La$ is that 
satisfies $A=A^\star$ where $A^\star$ is the transpose of $A$ followed by the 
natural involution $g \in F \to g^{-1} \in F$), defined up to a minor
indeterminacy that depends on the choice of lifts $\ti C_i$ of the components
of $C$. A different choice of lifts replaces the linking matrix 
$\lk(\ti C)$ by $P \, \lk(\ti C) P^\star$, where $P$ is a diagonal matrix with
elements of $\pm F \subset \La$ in the diagonal.

\begin{definition}
\lbl{def.W}
Taking into account the Kirby moves, we can define a map
$$
W: \NO/\la \kappa \ra \longto \B(\La\to\BZ)
$$
where $W(C)=\lk(\ti C)$ and where 
$\B(\La\to\BZ)$ denotes the set of {\em simple stable congruence classes
of Hermitian matrices over $\La$, invertible over $\BZ$}\footnote{By 
invertible over $\BZ$, we mean that the image of a matrix which is
defined over $\La$ under the map $\La\to\BZ$ is invertible over $\BZ$.}.
Following \cite{GL1}, we call two Hermitian matrices $A,B$ 
{\em simply stably congruent} iff 
$A\oplus S_1=P(B\oplus S_2)P^\star$ for some diagonal matrices $S_1,S_2$
with $\pm 1$ entries and some matrix $P$ which is either {\em elementary}
(i.e., it differs from the identity matrix on a single non-diagonal entry) 
or {\em diagonal} (with entries in $\pm F \subset \La$).
\end{definition}

Let $\Lhat$ denote the completion of the rational group-ring $\BQ[F]$
with respect to the augmentation ideal, that is the ideal generated by
$t_1-1$, where $t_1,\dots,t_g$ generate $F$. $\Lhat$ can be identified with 
the ring $\BQ[[h_1,\dots,h_g]]$ of formal power series in non-commuting
variables $h_i$. 

\begin{definition}
\lbl{def.chi}
Let 
$$ 
\chi: \Mat(\Lhat\to\BZ) \longto \Lhat 
$$ be the invariant on the set
$\Mat(\Lhat\to\BZ)$ 
{\em of matrices $W$ over $\Lhat$, nonsingular over 
$\BZ$}, defined by
$$\chi(W)=\tr\logp(W)
$$
where
$$\logp(W)=\sum_{n=1}^\infty\frac{(-1)^n}{n}\tr
\left((W (\e W)^{-1}-I)^n\right)
$$
and $\e:\La\to\BZ$. 
\end{definition}

In a sense, $\chi$ is a logarithmic determinant. In our case, we
want to define $\chi$ on the quotient $\Bla$ of the set of hermitian
matrices over $\La$ which are invertible over $\BZ$.
Unfortunately, $\chi$ does not descend to a $\Lhat$-valued invariant on
the set $\Bla$, but it does descend to a $\Lhat/(\mathrm{cyclic})$-valued
invariant, where $\Lhat/(\mathrm{cyclic})$ is the abelian group quotient
of the abelian group $\Lhat$ modulo the subgroup $xy-yx$ for $x,y \in
\Lhat$. Equivalently, $\Lhat/(\mathrm{cyclic})$ can be described as the
free abelian group generated by cyclic words in $h_i$ variables. 

Our next result shows that the $\chi$ invariant is a noncommutative analogue 
of the Alexander polynomial of an $F$-link.
As a motivation, recall that the Alexander polynomial of a knot measures the 
order of the first homology of the free abelian cover of its complement.
Given an $F$-link $(M,L,\phi)$, we can consider the free cover $\ti X$
of its complement $X=M\sminus L$, and its homology $H_1(\ti X)$, 
which is a module over the group-ring $\La$. Since the ring $\La$
is not commutative (when $L$ is not a knot), the notion of a torsion module 
does not make sense. However, in \cite{GL2}, Levine and the first author 
defined an invariant $\chiD$ of $F$-links with values in $\Lhat/(\cyclic)$, 
the abelian group of formal power series (with rational coefficients) 
in noncommutative variables, modulo cyclic permutation. The invariant $\chiD$ 
is a reformulation of an invariant defined earlier by Farber, \cite{Fa}, and
in the case of knots equals, up to a multiple, to the logarithm of the
Alexander polynomial. Thus, $\chiD$ can be thought as a {\em noncommutative
torsion} invariant of $F$-links.

The next result gives an independent definition of the noncommutative torsion
of an $F$-link, and of its Blanchfield form.

\begin{theorem}
\lbl{thm.classical}
{\rm (a)} There exist an invariant $W$ of $F$-links 
with values in $\B(\La\to\BZ)$. \newline
{\rm (b)} $H_1(\ti X)$ is a $\La$-module with presentation $W$.
\newline
{\rm (c)} The Blanchfield pairing \cite{Dv}
$$
H_1(\ti X,\BZ)\otimes H_1(\ti X,\BZ) \to \Lloc/\La
$$
is given by $W^{-1} \bmod \La $, where $\Lloc$ is the Cohn localization of 
$\La$, \cite{FV}. \newline
{\rm (d)} $\chi=\chi_{\D}$.  \newline
\end{theorem}

Thus, the $W$ invariant of an $F$-link determines its Blanchfield form.
Conversely, using various exact sequences of algebraic surgery, one can show 
that the Blanchfield form determines $W$, as was communicated to us by
D. Sheiham. Since this would take us far afield, we prefer to postpone this 
to a future publication.

\section{proofs}
\lbl{sec.proofs}

\subsection{Beginning the proof of Theorem \ref{thm.surgeryF}}
\lbl{sub.surgery}

Fix once and for all a {\em based unlink} $\O$ of $g$ components in $S^3$.
A based unlink consists of an unlink $\O$ in $S^3$, together with arcs
from a base point in $S^3-\O$ to each of the component of the unlink,
as shown below (for $g=3$):
$$
\psdraw{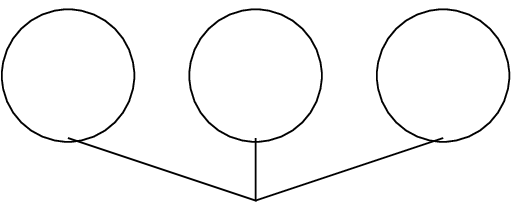}{1in}
$$ 
Let $\NO$ denote the set of {\em null homotopic links $L$ with 
$\BZ$-invertible linking matrix} in the complement of $\O$.

Surgery on an element $C$ of $\NO$ transforms $(S^3,\O)$ to a $F$-link $(M,L)$.
Indeed, since $C$ is null homotopic the natural map $\pi_1(S^3\sminus\O)\to F$
gives rise to a map $\pi_1(M\sminus L)\to F$. Alternatively, one can construct
disjoint Seifert surfaces for each component of $L$ by tubing the disjoint
disks that $\O$ bounds, which is 
possible, since each component of $C$ is null homotopic. Since the linking
matrix of $C$ is invertible over $\BZ$, $M$ is an \ihs .
Let $\stackrel{\kappa}\sim$ denote the equivalence relation on $\NO$ 
generated by the moves of {\em handle-slide} $\kappa_1$ 
(i.e., adding a parallel of a link component to another component) and 
{\em stabilization} $\kappa_2$ (i.e., 
adding to a link an unknot away from the link with framing $\pm 1$). 
It is well-known that $\kappa$-equivalence preserves surgery, thus
we get a surgery map
$$
\NO/\la \kappa \ra \longto F-\links .
$$ 

Now, we will show that the surgery map is onto. Given
a $F$-link $(M,L,\phi)$ choose disjoint Seifert surfaces $\S$ for its 
components
and choose a framed link $C$ in $M$ such that $M_{C}=S^3$. $C$ might
intersect $\S$, however by a small isotopy (which does not change the 
property that $M_{C}=S^3$) we may assume that $C$ is disjoint from $\S$.
Consider the image $\S'$ of $\S$ in $S^3=M_{C}$.
By reversing the $C$ surgery on $M$, consider the image $C'$ of a parallel
of $C$ in $S^3=M_C$ such that $M=S^3_{C'}$.   
By doing crossing changes among the bands of $\S'$, translated in terms of 
surgery on a unit-framed link $C'' \subset S^3\sminus\S'$, we can assume 
that each component of $\S'$ lies in a ball disjoint from the other 
components and further that the boundary of $\S'$ is an unlink. 
The link $C' \cup C''$ is null homotopic in the complement of the unlink 
$\O=\pt\S'$ since it is disjoint from $\S'$, and since the map 
$\phi:\pi_1(S^3\sminus\O)\cong F$ is given geometrically by intersecting
a loop in $S^3\sminus\O$ with $\S'$ and recording the corresponding element
of the free group. 

The linking matrix of 
$C' \cup C''$ is invertible over $\BZ$ since $S^3_{C' \cup C''}=M$,
an \ihs . Now perform an isotopy and choose arcs $\ga$ from a base point to
each of the components of $\S'$ in such a way that
$\S' \cup \ga$ intersects a ball as follows (for $g=3$): 
$$ 
\psdraw{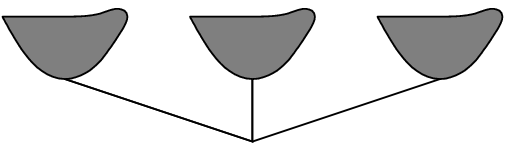}{1in}
$$
By construction, 
$C' \cup C'' \in \NO$ maps to $(M,L,\phi)$.

Next, we need to show that the surgery map is 1-1; we will do this in the
next section.

\subsection{A relative version of Fenn-Rourke's theorem}
\lbl{sub.relative}

The goal of this independent section is to show a  relative
version of Fenn-Rourke's theorem \cite[Theorem 6]{FR}. We begin with some
notation. Given a 3-manifold $M$ (possibly with nonempty boundary)
and a framed link $L$ in its interior, surgery produces a 3-manifold $M_L$
whose boundary is canonically identified with that of $M$
and a cobordism $W_L$ with boundary $M \cup (\partial M \times I) \cup M_L$.

\begin{theorem}
\lbl{thm.FR}
If M is a 3-manifold with boundary, $L_1$, $L_2$ two framed links 
in its interior, then $L_1$ and $L_2$ are $\kappa$-equivalent iff there exist
\begin{itemize}
\item[{\rm (a)}]
$h: M_{L_1}\to M_{L_2}$, a diffeo rel. boundary.
\item[{\rm (b)}]
An isomorphism $\iota: \pi_1(W_{L_1})\longto\pi_1(W_{L_2})$.
\item[{\rm (c)}]
A commutative diagram $(\Delta)$
$$
\divide\dgARROWLENGTH by2
\begin{diagram}
\node{\pi_1(M_{L_1})}       
\arrow[2]{e,t}{h_\ast}\arrow{s}
\node[2]{\pi_1(M_{L_2})}
\arrow{s}                   \\
\node{\pi_1(W_{L_1})}       
\arrow[2]{e,t,..}{\iota}
\node[2]{\pi_1(W_{L_2})}    \\
\node[2]{\pi_1(M)}
\arrow{nw}\arrow{ne}
\end{diagram}
$$
\end{itemize}
such that $\eta(\Delta)=0 \in H_4(\pi)$.
\end{theorem}
Here, $\eta(\D) \in H_4(\pi)$ is the homology class of $\rho(W)$, where
$\rho: W \to K(\pi,1)$, $W=W_{L_1} \cup_\partial W_{L_2}$ (where the 
boundaries are identified via the identity on $M$, $\partial M \times I$
and $h$ on the top piece) and 
$\pi=\pi_1(W_{L_1})= \pi_1(W_{L_2})$. Note that $W$ is a closed 
4-manifold. 

The proof of Theorem \ref{thm.FR} follows mainly from Fenn-Rourke \cite{FR}, 
with some simplifications due to Roberts, \cite{Rb}.  
We begin with the following lemma from \cite{FR}, which reformulates
the vanishing of the condition $\eta(\D)=0$.

\begin{lemma}\cite[Lemma 9]{FR}
\lbl{lem.FR}
Assume that we are in the situation of Theorem \ref{thm.FR} and that there
exists $h,\iota$ making $\D$ commute. Then $\eta(\D)=0$ iff 
$W\#_r \BC P^2 \#_s \overline{\BC P^2}=\partial \Omega^5$ for some
$r,s$ and closed 5-manifold $\Omega$ such that
the diagram
$$
\divide\dgARROWLENGTH by2
\begin{diagram}
\node{\pi_1(W_{L_1})}
\arrow[2]{e,t}{\iota}\arrow{se,r}{(j_1)_\ast}  
\node[2]{\pi_1(W_{L_2})} 
\arrow{sw,r}{(j_2)_\ast} \\
\node[2]{\pi_1(\Omega^5)}
\end{diagram}
$$
commutes and $(j_1)_\ast, (j_2)_\ast$ are split injections.
\end{lemma}

\begin{proof}(of Theorem \ref{thm.FR})
If $L_1$ and $L_2$ are $\kappa$-equivalent, then it follows from the ``only
if'' statement of Lemma \ref{lem.FR} that $\eta(\D)=0$ since the induced 
homeomorphisms on $\pi_1$ do not change.

So, assume that the algebraic conditions of Theorem \ref{thm.FR} are satisfied
for $L_1$, $L_2$. We will show that $L_1$ is $\kappa$-equivalent to $L_2$.

The first step is to arrange via $\kappa$-moves that $W_{L_1}$
and $W_{L_2}$ diffeomorphic. The proof of \cite[p.9]{FR} (where they
assume $\partial M= \phi$) using surgery works here, too.

The next step is to use Cerf theory. Namely, consider Morse functions $f_i$ 
realizing the cobordisms $W_{L_i}$ for $i=1,2$. Since $W_{L_1}$ is 
diffeomorphic to $W_{L_2}$, $f_1$ and $f_2$ are related by a one parameter
family of smooth functions $f_t$ for $1 \leq t \leq 2$. The arguments of 
Roberts \cite[Stage 2, p.3]{Rb} in this case imply that $L_1$ is 
$\kappa$-equivalent to $L_2$. 
\end{proof}

\begin{remark}
\lbl{rem.FR}
If $L$ is a null homotopic link in $M$, then $W_L$ is homotopic to a one point
union of $M$ and some number of 2-spheres $S^2$.
In particular, $\pi_1(M)\simeq \pi_1(W_L)$.
\end{remark}

{\em End of the Proof of Theorem \ref{thm.surgeryF}.} 
We will apply Theorem \ref{thm.FR} to $M=S^3\sminus\O$. 
Consider $(L_i,\phi_i)$ 
two surgery presentations of the same $F$-link for $i=1,2$. Then, we have a 
homotopy commutative diagram
$$
\divide\dgARROWLENGTH by2
\begin{diagram}
\node{M_{L_1}}    
\arrow{e,t}{h}\arrow{s,l}{\phi_1}
\node{M_{L_2}}  
\arrow{s,r}{\phi_2} \\
\node{M}          
\arrow{e,=}
\node{M}
\end{diagram}
$$
where $\phi_i:M_{L_i}\to M$ are splitting maps. It follows by Remark 
\ref{rem.FR} that a map $\iota$ exists such that the
following diagram $(\D)$ commutes
$$
\divide\dgARROWLENGTH by2
\begin{diagram}
\node{\pi_1(M_{L_1})}       
\arrow[4]{e,t}{h_\ast}\arrow{s}\arrow[2]{se,l,1}{\phi_1}
\node[4]{\pi_1(M_{L_2})}
\arrow{s}\arrow[2]{sw,l,1}{\phi_2}  \\
\node{\pi_1(W_{L_1})}       
\arrow{e,-,..}
\node{}
\arrow[2]{e,t,-,..}{\iota}
\node[2]{}
\arrow{e,..}
\node{\pi_1(W_{L_2})}    \\
\node[3]{\pi_1(M)}          
\arrow{nww,r}{\sim}\arrow{nee,r}{\sim}
\end{diagram}
$$
The obstruction $\eta(\D)$ lies in $H_4(F)$. However, 
$K(F,1)=\vee S^1$, thus $H_4(F)=H_4(\vee S^1)=0$.
Thus, Theorem \ref{thm.FR} implies that $L_1$ is $\kappa$-equivalent to
$L_2$.
\qed

\subsection{A matrix-valued invariant of $F$-links}
\lbl{sub.chi}

The goal of this section is to prove Theorem \ref{thm.classical}.
Combining the map $\NO/\la \kappa \ra
\longto \Bla$ of Definition \ref{def.W} with the surgery map of Theorem
\ref{thm.surgeryF}, gives rise to the map $W$ of Theorem \ref{thm.classical} 
(a).

Recall the $\chi$-invariant of Definition \ref{def.chi}.

\begin{proposition}
\lbl{prop.w} 
{\rm (a)} For matrices $A,B$ over $\Lhat$, nonsingular over $\BZ$, we
have in $\Lhat/(\mathrm{cyclic})$ that
$$\chi(AB)=\chi(A)+\chi(B) 
\,\,\,\, \text{ and } \,\,\,\, \chi(A\oplus B)=\chi(A)
+\chi(B).$$
{\rm (b)} $\chi$ descends to a $\Lhat/(\mathrm{cyclic})$-valued invariant
of the set $\B(\La\to\BZ)$. \newline
{\rm (c)} For $A$ as above, $\chi(A)=\chi(A \, \e (A^{-1}))$ where 
$\e:\La\to\BZ$. \newline
{\rm (d)} If $A$ is a matrix with integer entries
invertible over $\BZ$, then $\chi(A)=0$.
\end{proposition}

\begin{proof}
{\rm (b)} follows from {\rm(a)}. For \rm{(a)}, let $\R$ denote the ring of 
square matrices
with entries in $\Lhat$, and $[\R,\R]$ denote the subgroup of it generated
by matrices of the form $AB-BA$ where $A,B \in \R$.
Given two matrices $A,B \in \R$, with augmentations $A_0=\e A$ and
$B_0=\e B$, we claim that 
\begin{equation}
\lbl{eq.logdet}
\logp(AB)=\logp(A)+\logp(B) \,\, \bmod [\R,\R].
\end{equation}
Indeed, if $A_0=B_0=I$, then this follows from the Baker-Cambell-Hausdorff
formula, since $\logp(AB)=\log(AB)$. For the general case, we have
modulo $[\R,\R]$, that
\begin{eqnarray*}
\logp(AB) &=& \log(AB(A_0B_0)^{-1})=\log(ABB_0^{-1}A^{-1}AA_0^{-1}) \\
&=&
\log(ABB_0^{-1}A^{-1})+\log(AA_0^{-1}) \\
&=& A\log(BB_0^{-1})A^{-1}+
\log(AA_0^{-1})=\logp(B)+\logp(A).
\end{eqnarray*}
Since
$$
\tr(AB)=\tr(BA) \in \Lhat/(\cyclic),
$$
after taking traces, Equation \eqref{eq.logdet} 
implies {\rm (a)}. \newline
{\rm (c)} follows immediately from the definition of $\chi$ and the fact 
that $\e(A \, \e A^{-1})=I$. \newline
{\rm (d)} follows from that fact that if $A$ is invertible over $\BZ$,
then $A=\e A$, which implies $\logp(A)=0$, thus $\chi(A)=\tr\logp(A)=0$.
\end{proof}

We will also denote by $\chi$  the $\Lhat/(\mathrm{cyclic})$-valued
invariant of $F$-links defined by 
$$F-\links \cong \NO/\la\kappa\ra\to\Bla\to\Lhat/(\mathrm{cyclic}).$$

\begin{proof}(of Theorem \ref{thm.classical})
Consider $C \in \NO$ and let $(M,L,\phi)$ denote the corresponding $F$-link. 
Let $W$ denote the equivariant linking matrix of $C$, and $\ti X$ 
(resp. $\ti X_0$) denote the free cover of $X=M\sminus L$ (resp.
$X_0=S^3\sminus\O$) given by $\phi$. Note that $\ti X$ is obtained
from $\ti X_0$ by surgery on $\ti C$.

Part (b) follows from classical surgery arguments presented in the 
case of knots by Levine in \cite{Le1,Le2} and adapted without difficulty in 
the case of $F$-links; see also \cite[Lemma 2.4]{GL3}. For the case of knots,
Lemma 12.2 and the discussion of p.46-47 in \cite{Le2}, imply that the 
Blanchfield pairing is given by $W^{-1}$ modulo $\La$. This discussion
can be generalized without change to the case of $F$-links, and proves
part (c).

For part (d), observe that the $\chi$ and the $\chiD$ invariant of \cite{GL2}
are defined, respectively, via the Seifert surface and the surgery view of
boundary links. The two views can be related using the so-called Y-view of 
knots presented in \cite[Section 6.4]{GGP} and extended without difficulty to 
the case of boundary links. Consider a disjoint union $\S$ of Seifert surfaces
in $S^3$ which we think of as a disjoint union of embedded disks with pairs 
of bands attached in an alternating way along each disk. 

Consider an additional link $L'$ in $S^3\sminus \S$ such that its linking
matrix $C$ satisfies $\det(C)
=\pm 1$ and such that the linking number between the cores of the bands
and $L'$ vanishes. With respect to a choice of orientation 
of 1-cycles corresponding to the cores of the bands, a Seifert matrix of 
$\S$ is given by
$$
A=\mtwo \Lxx {\Lxy +I} \Lyx \Lyy
$$
Start from a surgery presentation of $(S^3,\S)$ in terms of {\em clovers with 
three leaves}, as was explained in \cite[Section 6.4]{GGP} and summarized
in the following figure:
$$ 
\psdraw{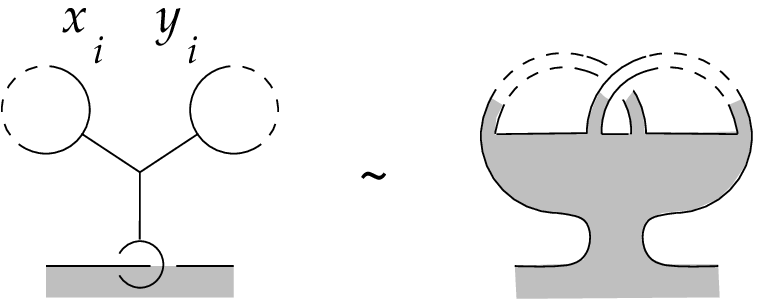}{2in} 
$$
Surgery on a clover with three leaves can be described in terms of surgery
on a six component link $L'''$. It was observed by the second author in 
\cite[Figure 3.1]{Kr} that $L'''$ can be simplified via Kirby moves to a four 
component link $L''$ disjoint from a surface $\S'$ shown 
as follows (where $\{x_i,y_i\}$ correspond to the same link in the
figure above and below):
$$ \psdraw{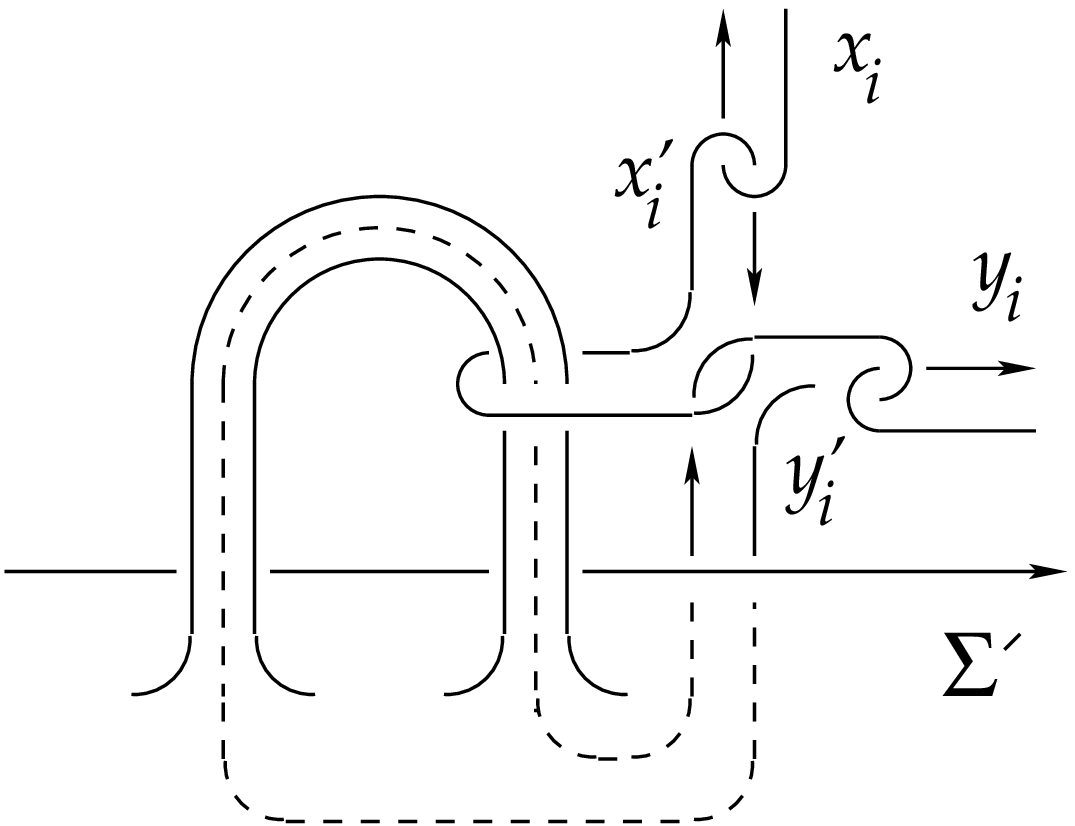}{1.7in} $$
The equivariant linking matrix of (a based representative of) $L'' \cup L'$ 
is given by $W \oplus C$ where $W$ is the following matrix 
in $\{x_i,y_i,x_i',y_i'\}$ basis:

$$
W=\begin{pmatrix} 
\Lxx & \Lxy & I & 0 \\
\Lxy & \Lyy & 0 & I \\
I    & 0    & 0 &{xI} \\
0 & I &  {\fx I} & 0
\end{pmatrix} 
= \mtwo L I I  B
$$
where $x=t-1$, $\fx=\ft -1$, $\ft=t^{-1}$ and
$$
B=\mtwo 0 {xI} {\fx I} 0 .
$$
Furthermore, 
$$
W W(1)^{-1} =
\mtwo L I I B \mtwo 0 I I {-L} = \mtwo I 0 B {I-BL} .
$$
thus, by part (c) of Proposition  $\chi(W)=\chi(I-BL)$. 

It is easy to see that every $F$-link has a surgery presentation $L'' \cup L'$
as above by adapting the onto part of the proof of Theorem \ref{thm.surgeryF}.
We can now compute as follows:
\begin{xalignat*}{2}
\chi(W \oplus C) &= \chi(W)+ \chi(C) & \text{by Proposition \ref{prop.w}(a)}
\\
&= \chi(W) & \text{by Proposition \ref{prop.w}(d), since $C$ is invertible}
\\
&= \chi(W W(1)^{-1}) & \text{by Proposition \ref{prop.w}(c)} 
\\
&= \chi(I-BL) & \text{by the above calculation} 
\end{xalignat*}

Using the notation of \cite{GL2}, let $Z=A(A-A')^{-1}$. Since 
$(A-A')^{-1}=\mtwo 0 {-I} I 0$   it follows that 
$$
L + \mtwo 0 I 0 0 = A = Z \mtwo 0 I {-I} 0 \,\, \text{ thus } \,\,
L=Z \mtwo 0 I {-I} 0 - \mtwo 0 I 0 0 .
$$
Substituting for $L$,  we obtain that
\begin{eqnarray*}
I-BL &=& \mtwo I 0 0 I -B\left( Z \mtwo 0 I {-I} 0 - \mtwo 0 I 0 0 \right)
\\
&=&
\mtwo I 0 0 {\ft I} - \mtwo 0 {x I} {\fx I} 0 Z \mtwo 0 I {-I} 0 \\
&=&
\left( \mtwo 0 {-I} {\ft I} 0 - \mtwo 0 {x T} {\fx I} 0 Z \right)
\mtwo 0 I {-I} 0 \\ 
&=&
\mtwo 0 {-I} {\ft I} 0 (I + XZ) \mtwo 0 I {-I} 0 
\end{eqnarray*}
where $X=\mtwo {x I} 0 0 {x I}$.
Proposition \ref{prop.w} implies that 
\begin{xalignat*}{2}
\chi(I-BL) &= \chi \mtwo 0 {-I} {\ft I} 0 + \chi (I+XZ)
+ \chi \mtwo 0 I {-I} 0 & \\
&= \chi (I+XZ) & \\
&= \chi_{\D}(A) & \text{by the definition of $\chi_{\D}$}
\end{xalignat*}
as required.
\end{proof}

\subsection{Acknowledgement}
The second author wishes to thank D. Bar-Natan for support and encouragement.

\ifx\undefined\bysame
	\newcommand{\bysame}{\leavevmode\hbox
to3em{\hrulefill}\,}
\fi

\end{document}